\documentclass[12pt,a4paper,psamsfonts]{amsart} 

\usepackage{srcltx}
\usepackage[T1]{fontenc}
\usepackage[latin1]{inputenc}
\usepackage{enumerate}
\usepackage{graphicx}
\usepackage{psfrag}
\usepackage{pst-all} 
\usepackage{pstricks-add}
\usepackage{setspace}
\usepackage{cite}

\usepackage[neverdecrease]{paralist} 
\newcommand\enumref[1]{\eqref{#1}} 

\usepackage[vmargin={30mm,30mm},hmargin={30mm,30mm}]{geometry}  

\newtheorem{theorem}{Theorem}[section]

\newtheorem{lemma}[theorem]{Lemma}
\newtheorem{proposition}[theorem]{Proposition}
\newtheorem{corollary}[theorem]{Corollary}
\theoremstyle{definition}

\newtheorem{example}{Example}
\theoremstyle{remark}
\newtheorem{remark}{Remark}

\newtheoremstyle{citing}
  {3pt}
  {3pt}
  {\itshape}
  {}
  {\bfseries}
  {.}
  {.5em}
  {\thmnote{#3}}
\theoremstyle{citing}

\allowdisplaybreaks[3] 

\DeclareMathOperator{\supp}{supp} %

\DeclareMathOperator{\exponential}{e}

\DeclareMathOperator{\determinant}{det}

\newcommand\Frame{S}
 \newcommand{\tran}[1][\gamma]{T_{#1}} 
 \newcommand{\dila}[1][A^j]{D_{\! #1}} 

\newcommand{\expo}[1]{\exponential^{#1}}

\newcommand{\dtm}{\ensuremath{\determinant}}

\newcommand*{\numbersys}[1]{\ensuremath{\mathbb{#1}}} 
\newcommand*{\C}{\numbersys{C}}
\newcommand*{\R}{\numbersys{R}}

\newcommand*{\Z}{\numbersys{Z}}
\newcommand*{\N}{\numbersys{N}}

\newcommand*{\Rn}{\numbersys{R}^n}

\newcommand*{\Zn}{\numbersys{Z}^n}
\newcommand{\itvoc}[2]{\ensuremath{\left({#1},{#2}\right]}} 
\newcommand{\itvoo}[2]{\ensuremath{\left({#1},{#2}\right)}} %
\newcommand{\itvcc}[2]{\ensuremath{\left[{#1},{#2}\right]}} %
\newcommand{\itvco}[2]{\ensuremath{\left[{#1},{#2}\right)}} %
\newcommand{\itvccs}[2]{\ensuremath{\lbrack{#1},{#2}\rbrack}} %
\newcommand{\abs}[1]{\ensuremath{\left\lvert#1\right\rvert}}
\newcommand{\absBig}[1]{\ensuremath{\Bigl\lvert#1\Bigr\rvert}}
\newcommand{\abssmall}[1]{\ensuremath{\lvert#1\rvert}}
\newcommand{\norm}[2][]{\ensuremath{\left\lVert#2\right\rVert_{#1}}}
\newcommand{\normsmall}[2][]{\ensuremath{\lVert#2\rVert_{#1}}}
 
\newcommand{\enorm}[2][]{\abs{#2}_{#1}} 
\newcommand{\enormsmall}[2][]{\abssmall{#2}_{#1}} 
\newcommand{\einnerprod}[3][]{\ensuremath{\left\langle
      #2,#3\right\rangle_{\! #1}}} 
\newcommand{\innerprod}[3][]{\ensuremath{\left\langle #2,#3\right\rangle_{\! #1}}}

\newcommand{\set}[1]{\ensuremath{\left\lbrace{#1}\right\rbrace}}
\newcommand{\setprop}[2]{\ensuremath{\left\lbrace{#1} : {#2}\right\rbrace}}

\newcommand{\setsmall}[1]{\ensuremath{\lbrace{#1}\rbrace}}

\newcommand{\setpropsmall}[2]{\ensuremath{\lbrace{#1} : {#2}\rbrace}}
\usepackage{xspace}

\newcommand{\ie}{\textit{i.e.,\xspace}}
 
\newcommand{\eg}{\textit{e.g.\@\xspace}}
\newcommand{\almoste}{a.e.} 

\newcommand{\lat}[1]{\ensuremath {#1}} 
\newcommand{\LG}{\ensuremath\lat{\Gamma}}
\newcommand{\LL}{\ensuremath\lat{\Lambda}}
 
\newcommand{\unitball}[2]{\ensuremath{B(#1,#2)}}
\newcommand{\nn}{n}

\newcommand{\dd}{d}
\newcommand{\kstar}{\gamma}

\newcommand{\tranl}[1][bk]{T_{#1}}
 \renewcommand{\tran}[1][\gamma]{T_{#1}} 
 \renewcommand{\dila}[1][A^j]{D_{\! #1}} 

\usepackage[hyperindex,colorlinks,citecolor=blue]{hyperref} 
\usepackage{color}
\definecolor{darkblue}{rgb}{0,0,.5}
\definecolor{darkred}{rgb}{.3,0,0}
\hypersetup{colorlinks=true, breaklinks=true,
  linkcolor=darkred, menucolor=darkred, pagecolor=darkred,
  urlcolor=darkred, citecolor=darkblue, filecolor=darkred,  bookmarks=true,
bookmarksopen=true, bookmarksopenlevel=3, plainpages=false,
pdfpagelabels=true}

\hypersetup{
pdfview={FitH}, 
pdfstartview={FitH},
pdfauthor = {Jakob Lemvig},
pdftitle = {Constructing pairs of dual bandlimited  frame wavelets   
 in L^2(R^n) Wavelet frames and their duals},
pdfsubject = {Wavelet frames},
pdfkeywords = {real
   and expansive dilation matrix,  bandlimited wavelets, dual frames, non-tight
   frames, partition of unity},
pdfcreator = {LaTeX with hyperref package},
pdfproducer = {dvips + ps2pdf}}

\title[Constructing pairs of 
  dual bandlimited  frame wavelets   
 in $L^2(\mathbb{R}^n)$]
{Constructing pairs of 
  dual bandlimited \\  frame wavelets   
 in $L^2(\mathbb{R}^n)$}

\author{Jakob Lemvig}
\address{Institute of Mathematics, University of Osnabr\"uck,
    49069 Osnabr\"uck, Germany}
    \email{jlemvig@uni-osnabrueck.de}

  \date{\today}

\begin{document}

\begin{abstract}
  Given a real, expansive dilation matrix we prove that any
  band\-li\-mi\-ted function $\psi \in L^2(\Rn)$, for which the
  dilations of its Fourier transform form a partition of unity,
  generates a wavelet frame for certain translation lattices.
  Moreover, there exists  a dual wavelet frame generated by a finite linear
  combination of dilations of $\psi$ with explicitly given
  coefficients. The result allows a simple construction procedure for
  pairs of dual wavelet frames whose generators have compact support
  in the Fourier domain and desired time localization. The
  construction relies on a technical condition on $\psi$, and we
  exhibit a general class of  function satisfying this condition. 
\end{abstract}

\maketitle

\section{Introduction}
\label{sec:intro}
For $A \in GL_\nn(\R)$ and  $y \in \Rn$, we define the dilation operator
on $L^2(\Rn)$ by
$\dila[A] f(x)=\abs{\dtm{A}}^{1/2}f(Ax)$ and the translation operator by
$\tran[y]f(x)=f(x-y)$. Given a $\nn\times \nn$ real, expansive
matrix $A$ and a  lattice of the form $\LG=P \Zn$ for
$P \in GL_\nn(\R)$, we consider wavelet systems of the form 
 \[ \set{\dila \tran \psi}_{j \in \Z, \gamma \in \LG}, \]
where the Fourier transform of $\psi$ has compact support. 
Our aim is, for any given real, expansive dilation matrix $A$, to construct
wavelet frames with good regularity properties and with a dual frame
generator of the form 
\begin{equation}
  \label{eq:26}
\phi = \sum_{j = a}^b c_j \dila \psi
\end{equation}
for some explicitly given coefficients $c_j \in \C$ and $a,b \in \Z$.
This will generalize and extend the one-dimensional results on
constructions of dual wavelet frames in \cite{lemvig_constr_paris,MR2311859} to
higher dimensions. The extension is non-trivial since it is unclear
how to determine the translation lattice $\LG$ and how to
control the support of the generators in the Fourier domain. This will
be done by considering suitable norms in $\Rn$ and non-overlapping
packing of ellipsoids in lattice arrangements. 

The construction of redundant wavelet representations in higher
dimensions is usually based on extension principles 
\cite{ehler_han,id_bh_ar_zs_03,MR1968120,MR1992289,MR2000c:42037,MR98c:42035,MR2355010,MR2274841,MR1848710}.  
By making use of extension principles one is restricted to considering
expansive dilations $A$ with integer coefficients. Our constructions
work for any real, expansive dilation. Moreover, in the extension
principle the number of generators often increases with the smoothness
of the generators. We will construct pairs of dual wavelet frames
generated by one smooth function with good time localization.

It is a well-known fact that a wavelet frame need not have dual frames
with wavelet structure. In \cite{MR2286929} frame wavelets with
compact support and explicit analytic form are constructed for real
dilation matrices. However, no dual frames are presented for these
wavelet frames. This can potentially be a problem because it might be
difficult or even impossible to find a dual frame with wavelet
structure. Since we exhibit \emph{pairs} of dual wavelet frames, this issue
is avoided.

The principal importance of having a dual
generator of the form (\ref{eq:26}) is that it will inherit properties
from $\psi$ preserved by dilation and linearity, \eg vanishing
moments, good time localization and regularity properties. For a more
complete account of such matters we refer to
\cite{lemvig_constr_paris}.

In the rest of this introduction we review basic definitions. A frame for a separable
Hilbert space $\mathcal{H}$ is a countable collection of vectors $\{f_j\}_{j \in
  \mathbb{J}}$ for which there are
constants $0 < C_1 \leq C_2 < \infty$ such that
\[ 
C_1 \norm{f}^2 \leq \sum_{j \in \mathbb{J}} \abs{\innerprod{f}{f_j}}^2 \leq
C_2 \norm{f}^2 \qquad\text{for all }f\in \mathcal{H}.
\]
If the upper bound holds in the above inequality, then $\{f_j\}$ is
said to be a Bessel sequence with Bessel constant $C_2$. 
For a Bessel sequence $\{f_j\}$ we define the frame operator by 
\[
\Frame\colon \mathcal{H} \to \mathcal{H}, \qquad Sf = \sum_{j \in
  \mathbb{J}} \innerprod{f}{f_j} f_j. 
\] 
This operator is bounded, invertible, and positive. A frame $\{f_j\}$
is said to be \emph{tight} if we can choose $C_1 = C_2$; this is
equivalent to $\Frame = C_1 I$ where $I$ is the identity operator.
Two Bessel sequences $\{f_j\}$ and $\{g_j\}$ are said to be
\emph{dual} frames if
\[ f = \sum_{j \in \mathbb{J}} \innerprod{f}{g_j}f_j \quad \forall f
\in \mathcal{H}. \] It can be shown that two such Bessel sequences are
indeed frames. Given a frame $\{f_j\}$, at least one dual always
exists; it is called the canonical dual and is given by $\{S^{-1} f_j
\}$. Only a frame, which is not  a basis, has several duals.

 For $f \in L^1(\Rn)$ the Fourier transform
is defined by $\hat f(\xi) = \int_{\Rn} f(x)\expo{-2 \pi i
  \einnerprod{\xi}{x}} \mathrm{d}x$ with
the usual extension to $L^2(\Rn)$. 

Sets in $\Rn$ are, in general, considered equal if they are equal 
up to sets of measure zero. 
The boundary of a set $E$ is denoted by $\partial E$, the interior by
$E^\circ$, and the closure by $\overline{E}$. Let $B\in GL_\nn(\R)$. A
\emph{multiplicative tiling set} $E$ for $\setpropsmall{B^j}{j\in \Z}$
is a subset of positive measure such that
\begin{gather}
  \absBig{\Rn \setminus \bigcup_{j \in \Z} B^j(E)} = 0 \quad \text{and} \quad
  \abs{B^j(E) \cap B^l(E)} =0 \quad \text{for $l\neq j$.}
\end{gather}
In this case we say that $\setprop{B^j(E)}{j \in \Z}$ is an almost
everywhere partition of $\Rn$, or that it tiles $\Rn$. A multiplicative
tiling set $E$ is \emph{bounded} if $E$ is a bounded set and $0 \notin
\overline{E}$. By $B$-dilative periodicity of a function $f \colon \Rn
\to \C$ we understand $f(x)=f(Bx)$ for \almoste\ $x \in \Rn$, and by a
$B$-dilative partition of unity we understand $\sum_{j \in \Z} f(B^j
x) =1$; note that the functions in the ``partition of unity'' are not
assumed to be non-negative, but can take any real or complex value.

A (full-rank) lattice $\LG$ in $\Rn$ is a point set of the form $\LG=P \Zn$ for
some $P \in GL_\nn(\R)$. The determinant of $\LG$ is
$d(\LG)=\abs{\dtm{P}}$; note that the generating matrix $P$ is not
unique, and that $d(\LG)$ is independent of the particular choice of $P$.

\section{The general form of the construction procedure}
\label{sec:general-form-Rn}
Fix the dimension $n \in \N$. We let $A\in GL_\nn(\R)$ be
expansive, \ie all eigenvalues of $A$ have absolute
value greater than one, and denote the transpose matrix by $B = A^t$.   
For any such 
dilation $A$, we want to construct a pair of
functions that generate dual wavelet frames for some translation
lattice.
Our construction is based on the following result which is a
consequence of the characterizing equations for dual wavelet frames by
Chui, Czaja, Maggioni, and Weiss \cite[Theorem 4]{MR1891728}.
\begin{theorem} \label{thm:dual-charac-Rn}
Let $A \in GL_\nn(\R)$ be expansive, let $\LG$ be a lattice in $\Rn$, and
let $\Psi=\{\psi_1, \dots, \psi_L\}$, $\tilde \Psi=\{\tilde\psi_1, \dots,
\tilde\psi_L\} \subset L^2(\Rn)$. Suppose that the two wavelet systems 
$\setpropsmall{\dila \tran \psi_l}{j\in \Z, \gamma \in \LG, l = 1, \dots, L}$
and $ \setpropsmall{\dila \tran \tilde{\psi}_l }{j\in \Z, \gamma\in \LG, l =
  1, \dots, L}$ form Bessel families. Then $\setsmall{\dila \tran
  \psi_l}$ and $\setsmall{\dila \tran \tilde\psi_l}$ will be dual
frames if the following conditions hold
   \begin{align}
    \label{eq:diagonalterm-Rn}
    &\sum_{l=1}^L \sum_{j\in \Z} \hat{\tilde{\psi_l}}(B^j \xi)
    \overline{\hat{\psi_l}(B^j \xi)} = d(\LG) &\quad &\text{\almoste\ } \xi
    \in
    \Rn, \\
\label{eq:nondiagonalterm-Rn}
    &\sum_{l=1}^L \hat{\tilde{\psi_l}}(\xi) \overline{\hat 
      \psi_l(\xi+\gamma)} = 0 &\quad &\text{\almoste\ } \xi \in
    \Rn \text{ for } \gamma \in \LG^\ast \setminus \{0\}.
  \end{align}
\end{theorem}
\begin{proof}
  By $\xi = B^j \omega$ for $j \in \Z$, condition
  (\ref{eq:nondiagonalterm-Rn}) becomes
\begin{align}\label{eq:1}
    \sum_{l=1}^L \hat{\tilde{\psi_l}}(B^j \omega) \overline{\hat 
      \psi_l(B^j \omega+\gamma)} = 0 \quad \text{\almoste\ } \omega \in
    \Rn \text{ for } \gamma \in \LG^\ast \setminus \{0\}.
  \end{align}
We use the notation as in \cite{MR1891728}, thus
$\Lambda (A, \LG) = \setpropsmall{\alpha \in \Rn}{ \exists (j, \gamma) \in \Z \times
  \LG^\ast :  \alpha = B^{-j}\gamma}$ and
$I_{A,\LG}(\alpha) = \setpropsmall{(j, \gamma) \in \Z \times \LG^\ast }{\alpha = B^{-j}\gamma}$.
Since $I_{A,\LG}(\alpha) \subset \Z \times (\LG^\ast \setminus \{0\})$ for any
$\alpha \in \Lambda(A,\LG) \setminus \{0\}$, equation~(\ref{eq:1}) yields
\begin{align*}
  \frac{1}{d(\LG)} \sum_{(j,\gamma)\in I_{A,\LG}(\alpha)}  \sum_{l=1}^L
  \hat{\tilde{\psi_l}}(B^j \omega) \overline{\hat  
      \psi_l(B^j( \omega+B^{-j}\gamma))} = 0 \quad \text{\almoste\ } \omega \in
    \Rn 
  \end{align*}
for $\alpha \neq 0$. By $I_{A,\LG}(0)= \Z \times \{0\}$, we can 
rewrite \eqref{eq:diagonalterm-Rn} as  
\begin{align*}
  \frac{1}{d(\LG)} \sum_{(j,\gamma)\in I_{A,\LG}(0)}  \sum_{l=1}^L
  \hat{\tilde{\psi_l}}(B^j \omega) \overline{\hat  
      \psi_l(B^j(\omega+B^{-j}\gamma))} = 1 \quad \text{\almoste\ } \omega \in
    \Rn,
  \end{align*}
using that $B^{-j}\gamma=0$ for all $j\in \Z$.
Gathering the  two equations displayed above yields
\begin{align*}
  \frac{1}{d(\LG)} \sum_{(j,\gamma)\in I_{A,\LG}(\alpha)}  \sum_{l=1}^L
  \hat{\tilde{\psi_l}}(B^j \omega) \overline{\hat  
      \psi_l(B^j(\omega+B^{-j}\gamma))} = \delta_{\alpha,0} \quad \text{\almoste\ } \omega \in
    \Rn,
  \end{align*}
for all $\alpha \in \Lambda(A,\LG)$. The conclusion follows now from
\cite[Theorem 4]{MR1891728}.
\end{proof}

The following result, Lemma~\ref{thm:bessel2-Rn}, gives a sufficient
condition for a wavelet system to form a Bessel sequence; it is an
extension of \cite[Theorem 11.2.3]{oc_03} from $L^2(\R)$ to
$L^2(\Rn)$.

\begin{lemma} \label{thm:bessel2-Rn} Let $A \in GL_\nn(\R)$ be 
  expansive, $\LG$ a lattice in $\Rn$, and $\phi \in L^2(\Rn)$.
  Suppose that, for some set $M \subset \Rn$ satisfying  $\cup_{l \in \Z} B^{l}(M) =
\Rn$,
\begin{align}\label{eq:14}
C_2 = \frac{1}{d(\LG)} \sup_{\xi \in M}  \sum_{j \in \Z}
\sum_{\kstar \in \LG^\ast}
\abs{\hat \phi(B^j \xi) \hat \phi(B^j \xi + \kstar)} < \infty{}.
\end{align}
Then the wavelet system  $\{\dila \tran \phi\}_{j \in \Z, \gamma \in
  \LG}$ is a Bessel sequence with bound $C_2$. Further, if also 
\begin{align}\label{eq:lemma-suff-lower}
C_1 = \frac{1}{d(\LG)} \inf_{\xi \in M} \left(\sum_{j \in \Z}
\abs{\hat \phi(B^j \xi)}^2  -  \sum_{j \in \Z}
\sum_{\kstar \in \LG^\ast \setminus \{0\}}
\abs{\hat \phi(B^j \xi) \hat \phi(B^j \xi + \kstar)} \right) > 0 ,
\end{align}
holds, then $\{\dila \tran \phi\}_{j \in \Z, \gamma \in \LG}$ is
a frame for $L^2(\Rn)$ with frame bounds $C_1$ and $C_2$.
\end{lemma}
\begin{proof}
The statement follows directly by applying Theorem 3.1 in \cite{chris_rahimi} on
generalized shift invariant systems to wavelet systems. 
In the general result for generalized
shift invariant systems \cite[Theorem 3.1]{chris_rahimi}, the supremum/infimum
is taken over $\Rn$, but  because of the $B$-dilative periodicity of
the series in (\ref{eq:14}) and (\ref{eq:lemma-suff-lower}) for wavelet systems, it 
suffices to take the supremum/infimum over a set $M \subset \Rn$ that has the
property that $\cup_{l \in \Z} B^{l}(M) = \Rn$ up to sets of measure
zero.
\end{proof}

Theorem~\ref{thm:dual-charac-Rn} and Lemma~\ref{thm:bessel2-Rn} are all we
need to prove the following result on pairs of dual wavelet frames.
\begin{theorem}\label{thm:constr-dual-wavelet-Rn-general} 
Let $A \in GL_\nn(\R)$ be expansive and  $\psi \in L^2(\Rn)$.
  Suppose that $\hat \psi$ is a bounded, real-valued function with $\supp \hat
  \psi \subset \cup_{j=0}^\dd B^{-j}(E) $ for some $d \in \N_0$ and
  some bounded multiplicative tiling set $E$ for  $\setprop{B^j}{j \in
    \Z}$, and  that  
  \begin{equation} 
   \label{eq:dyadic-part2-Rn}
   \sum_{j \in \mathbb{Z}} \hat \psi (B^j \xi) =1 
\quad \text{for \almoste\ }\xi \in \Rn.
\end{equation}
Let $b_j \in \C$ for $j=-\dd, \dots, \dd$ and let $\overline{m} =
\max\setprop{j}{b_j \neq 0}$ and $\underline{m} = -
\min\setprop{j}{b_j \neq 0}$. Take a lattice $\LG$ in $\Rn$ such that 
\begin{equation}
  \label{eq:13}
  \Bigl(\bigcup_{j=0}^\dd B^{-j}(E)+\kstar\Bigr) \cap \bigcup_{j=-\underline{m}}^{\overline{m}+\dd}B^{-j}(E) = \emptyset \quad
  \text{for all $\kstar \in \LG^\ast \setminus \{0\}$,}  
 \end{equation}
and define the function $\phi$ by
\begin{align}
  \label{eq:dual-generator1-Rn}
  \phi(x) = d(\LG) \sum_{j=-\underline{m}}^{\overline{m}} b_j \abs{\dtm{A}}^{-j} \psi(A^{-j}x) \quad
  \text{for } x \in \Rn.
\end{align}
If $b_0 = 1$ and $b_j + b_{-j} = 2$ for $j=1,2, \dots, \dd$, then
the functions $\psi$ and $\phi$ generate dual frames $\setsmall{\dila \tran \psi}_{j
  \in \Z, \gamma \in \LG}$ and $\setsmall{\dila \tran \phi}_{j \in \Z, \gamma
  \in \LG}$ for $L^2(\Rn)$.
  \end{theorem}

  \begin{proof}
 On the Fourier side, the definition in  (\ref{eq:dual-generator1-Rn}) becomes
\begin{align*}
  \hat \phi(\xi) = d(\LG) \sum_{j=-\underline{m}}^{\overline{m}} b_j \hat \psi(B^{j}\xi).
\end{align*}
Since $\hat \psi$ by assumption is compactly
supported in a ``ringlike'' structure bounded away from the origin,
this will also be the case for $\hat \phi$. This property
implies that $\psi$ and $\phi$ will generate wavelet Bessel sequences.
The details are as follows. The support of $\hat \psi$ and $\hat \phi$
is 
\begin{align}\label{eq:supp-generators-general}
  \supp \hat \psi \subset \bigcup_{j=0}^\dd B^{-j}(E), 
  &&\supp \hat \phi \subset
  \bigcup_{j=-\underline{m}}^{\overline{m}+\dd}B^{-j}(E).
\end{align}
Note that $0 \le  \underline{m}, \overline{m} \le
  \dd$.
The sets $\setprop{B^j(E)}{j \in \Z}$ tiles $\Rn$, whereby we see that
\begin{align}\label{eq:18}
\abs{\,\supp \hat \psi(B^j \cdot) \cap B^{-\dd}(E) }&= 0 \quad \text{for $j <0$
  and $j > \dd$,} \\
\intertext{and,} 
\abs{\,\supp \hat \phi(B^j \cdot) \cap B^{-\dd}(E) }&= 0 \quad \text{for $j
  <-\underline{m}$ and $j > \overline{m}+\dd$.} \label{eq:28}
\end{align} 
Since $\underline{m}, \overline{m} \ge 0$, condition (\ref{eq:13})
implies that $\hat \psi(B^j \xi) \hat \psi(B^j \xi + \kstar) = 0$ for
$j \ge 0$ and $\kstar \in \LG^\ast \setminus \{0\}$. Therefore, using 
 \eqref{eq:18}, we find that 
\begin{align*}
  \sum_{j\in\Z} \sum_{\kstar \in \LG^\ast } \abs{ \hat \psi(B^j \xi)
    \hat \psi(B^j \xi + \kstar)}
  =   \sum_{j=0}^{\dd} \left( \hat \psi(B^j \xi) \right)^2 < \infty
  \qquad \text{for $\xi \in B^{-\dd}(E)$.}
\end{align*}
An application of Lemma~\ref{thm:bessel2-Rn} with $M=B^{-\dd}(E)$ shows that
$\psi$ generates a Bessel sequence. Similar calculations using 
 \eqref{eq:28} will show that $\phi$ generates a Bessel sequence; in
this case the sum over $\kstar \in \LG^\ast$ will be finite, but it will
in general have more than one nonzero term.

To conclude that $\psi$ and $\phi$ generate dual wavelet frames we
will show that conditions~\eqref{eq:diagonalterm-Rn} and
\eqref{eq:nondiagonalterm-Rn} in Theorem~\ref{thm:dual-charac-Rn}
hold. By $B$-dilation periodicity of the sum in condition
\eqref{eq:diagonalterm-Rn}, it is sufficient to verify this condition
on $B^{-\dd}(E)$.
For $\xi  \in B^{-\dd}(E) $  we have by \eqref{eq:18},
\begin{align*}
\frac{1}{d(\LG)}  \sum_{j \in \Z} \overline{\hat \psi (B^j \xi)} \hat \phi (B^j \xi)
  & = \frac{1}{d(\LG)} \sum_{j =0}^{\dd} \hat \psi (B^j \xi) \hat \phi (B^j \xi) \\ &=
  \hat \psi (\xi) \left[ b_0 \hat \psi (\xi) + b_{1} \hat \psi(B \xi)
    + \dots + b_{\dd} \hat
    \psi(B^{\dd}\xi) \right] \\
  &\phantom{= }\;+ \hat \psi(B \xi) \left[ b_{-1} \hat \psi( \xi) +
    b_0 \hat \psi (B \xi) +
    \dots + b_{\dd-1}\hat  \psi(B^{\dd}\xi)  \right] + \cdots  \\
  &\phantom{= }\;+ \hat \psi(B^{\dd}\xi) \left[ b_{-\dd}\hat
    \psi(\xi) + \dots + b_{-1} \hat \psi(B^{\dd-1}\xi) + b_0 \hat
    \psi(B^{\dd}\xi) \right], 
\intertext{and further, by an expansion of these terms,  }
&=  \sum_{ j,l = 0}^{\dd} b_{l-j} \hat\psi(B^j \xi) \hat\psi(B^l \xi)\\
  &= b_0 \sum_{j=0}^{\dd} \hat\psi(B^j \xi)^2 + \sum_{ \substack{j,l = 0 \\ j >
  l}}^{\dd} (b_{j-l}+b_{l-j}) \hat\psi(B^j \xi) \hat\psi(B^l \xi) .
\end{align*}
Using that $b_0 = 1$ and $b_{j-l}+b_{l-j}= 2$ for $j\neq l$ and $j,l
=0,\dots, \dd$,  we arrive at
\begin{align*}
  \frac{1}{d(\LG)} \sum_{j \in \Z} \overline{\hat \psi (B^j \xi)} \hat
  \phi (B^j \xi) &= \sum_{j=0}^{\dd} \hat\psi(B^j \xi)^2 + \sum_{
    \substack{j,l = 0 \\ j >
      l}}^{\dd} 2\hat\psi(B^j \xi) \hat\psi(B^l \xi) \\
  &= \biggl(\sum_{j =0}^{\dd} \hat \psi (B^j \xi) \biggr)^2 =
  \biggl(\sum_{j \in \Z} \hat \psi (B^j \xi) \biggr)^2 = 1, 
\end{align*}
exhibiting that $\psi$ and $\phi$ satisfy condition \eqref{eq:diagonalterm-Rn}.

By \eqref{eq:supp-generators-general} we see that  condition
\eqref{eq:13} implies that  the functions $\hat
\phi$ and $\hat \psi (\cdot + \kstar)$ will have disjoint support
for $\kstar \in \LG^\ast \setminus \{0 \}$,
hence (\ref{eq:nondiagonalterm-Rn}) is satisfied.
\end{proof}

\begin{remark}
  The use of the parameters $b_j$ in the definition of the dual
  generator together with the condition $b_{-j}+b_j=2$ was first seen in
  the work of Christensen and Kim~\cite{OC_RYK_dual-Gabor-poly} on
  pairs of dual Gabor frames.
\end{remark}

We can restate
Theorem~\ref{thm:constr-dual-wavelet-Rn-general} for wavelet
systems with standard translation lattice $\Zn$ and dilation
$\widetilde A =P^{-1}AP$, where $P \in GL_n(\R)$ is so that $\LG = P
\Zn$. The result follows directly by an application of the relations
$\dila[\widetilde{A}^j] \dila[P] = \dila[P] \dila$ for $j \in \Z$
and $\dila[P]\tran[Pk]=\tran[k] \dila[P]$ for $k \in \Zn$, and the
fact that $\dila[P]$ is unitary as an operator on $L^2(\Rn)$.
\begin{corollary}
  Suppose $\psi$, $\{b_j\}$, $A$ and $\LG$ are as in
  Theorem~\ref{thm:constr-dual-wavelet-Rn-general}. Let $P \in
  GL_\nn(\R)$ be such that $\LG = P \Zn$, and let $\widetilde A
  =P^{-1}AP$. Then the functions $\tilde \psi = \dila[P]{\psi}$ and
  $\tilde \phi = \dila[P] \phi$, where $\phi$ is defined in
  (\ref{eq:dual-generator1-Rn}), generate dual frames
  $\setsmall{\dila[\widetilde A^j] \tran[k] \tilde \psi}_{j \in \Z, k
    \in \Zn}$ and $\setsmall{\dila[\widetilde{A}^j] \tran[k] \tilde
    \phi}_{j \in \Z, k \in \Zn}$ for $L^2(\Rn)$.
\end{corollary}

The following Example~\ref{ex:adhoc-construc-Rn} is an application of
Theorem~\ref{thm:constr-dual-wavelet-Rn-general} in $L^2(\R^2)$ for
the quincunx matrix. In particular, we construct a partition of unity
of the form (\ref{eq:dyadic-part2-Rn}) for the quincunx matrix. 
\begin{example}\label{ex:adhoc-construc-Rn}
  The quincunx matrix is defined as 
  \begin{equation*}
    A =
    \begin{pmatrix}
      1 & -1 \\ 1 & 1
    \end{pmatrix},
  \end{equation*}
  and its action on $\R^2$ corresponds to a counter clockwise rotation
  of 45 degrees and a dilation by $\sqrt{2} I_{2 \times 2}$.
\begin{figure}[htbp]
\psset{xunit=5cm,yunit=5cm}%
\begin{pspicture}(-.13,-.13)(1.2,1.2)
\psaxes[Dx=0.25, Dy=0.25,tickstyle=bottom]{->}(0,0)(-.13,-0.13)(1.2,1.2) 
\pspolygon(1,0)(1,1)(0.5,0.5)
\pspolygon(0,1)(1,1)(0.5,0.5)
\pspolygon(0,1)(0,0.5)(0.5,0.5)
\pspolygon(0.5,0)(0,0.5)(0.5,0.5)
\pspolygon(0.5,0)(1,0)(0.5,0.5)
\rput[c](0.5, 0.75){$J_5$} 
\rput[c]( 0.75, 0.5){$J_4$} 
\rput[c]( 0.625, 0.175){$J_2$} 
\rput[c]( 0.175, 0.625){$J_3$} 
\rput[c]( 0.35, 0.35){$J_1$} 
\uput[-90](1.17,-2.4pt){$x_1$}
\uput[145](-2.4pt,1.17){$x_2$}
\end{pspicture}\centering
\caption{Sketch of the triangular domains $J_i$, $i = 1,2,3,4,5$.} 
\label{fig:quincunx-support-sets}
\end{figure}
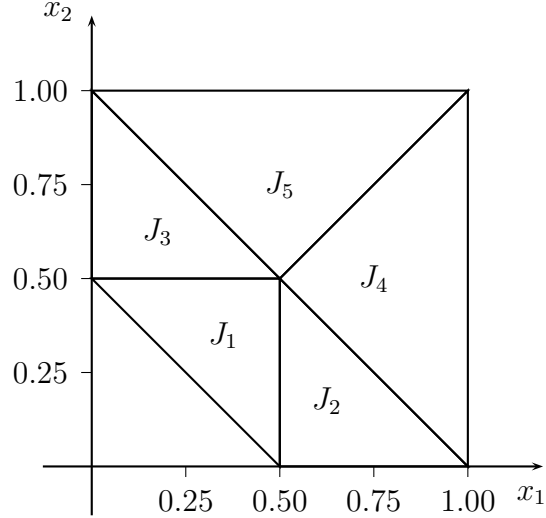
  Define the tent shaped,  piecewise linear function $g$ by 
  \begin{equation*}
    g(x_1,x_2)=\begin{cases}
      -1+2x_1 + 2x_2, \quad &\text{for $(x_1,x_2) \in J_1,$}\\
      2x_2, \quad &\text{for $(x_1,x_2) \in J_2,$}\\
      2x_1, \quad &\text{for $(x_1,x_2) \in J_3,$}\\
      2-2x_1, \quad &\text{for $(x_1,x_2) \in J_4,$}\\
      2-2x_2, \quad &\text{for $(x_1,x_2) \in J_5,$}\\
      0 & \text{otherwise,}
    \end{cases}
  \end{equation*}
where the sets $J_i$ are the triangular domains sketched in
Figure~\ref{fig:quincunx-support-sets}.   
Note that the value at ``the top of the tent'' is $g(1/2,1/2)=1$. 
Define $\hat\psi $ as a mirroring of $g$ in the $x_1$ axis and the $x_2$ axis:
\begin{equation*}
  \hat \psi(\xi_1, \xi_2) =
  \begin{cases}
    g(\xi_1,\xi_2)\quad &\text{for $(\xi_1,\xi_2) \in \itvco{0}{\infty}
      \times \itvco{0}{\infty},$}\\
g(\xi_1,-\xi_2)\quad &\text{for $(\xi_1,\xi_2) \in \itvco{0}{\infty}
      \times \itvoo{-\infty}{0},$}\\
g(-\xi_1,\xi_2)\quad &\text{for $(\xi_1,\xi_2) \in \itvoo{-\infty}{0}
      \times \itvco{0}{\infty},$}\\
g(-\xi_1,-\xi_2)\quad &\text{for $(\xi_1,\xi_2) \in \itvoo{-\infty}{0}
      \times \itvoo{-\infty}{0}.$}
  \end{cases}
\end{equation*}
Since the transpose $B$ of the quincunx matrix also corresponds to a
rotation of 45 degrees (but clockwise) and a dilation by $\sqrt{2} I_{2 \times 2}$, we
see that $ \sum_{j \in \mathbb{Z}} \hat \psi (B^j \xi) =1 $. 

We are now ready to apply
Theorem~\ref{thm:constr-dual-wavelet-Rn-general} with $E =
\itvcc{-1}{1}^2 \setminus B^{-1}(\itvcc{-1}{1}^2) = \itvcc{-1}{1}^2
\setminus I_1$ and $d=2$; the set $E$ is the union of the domians
$J_4$ and $J_5$ and their mirrored versions. We
choose $b_{-2}=b_{-1}=0$ and $b_1=b_2=2 d(\LG)$, hence
$\underline{m}=0$ and $\overline{m}=2$. Therefore, 
 \[ \bigcup_{j=0}^\dd B^{-j}(E),
 \bigcup_{j=-\underline{m}}^{\overline{m}+\dd}B^{-j}(E) \subset
 \itvcc{-1}{1}^2,\]
that shows that we can take 
$\LG^\ast = 2\Z^2$ or $\LG = 1/2 \Z^2$, since
$ (\itvcc{-1}{1}^2+
\kstar) \cap \itvcc{-1}{1}^2 = \emptyset$ whenever $0\neq \kstar \in
2\Z^2$. 
Defining the dual generator according to \eqref{eq:dual-generator1}
yields
\begin{equation}
  \label{eq:11}
  \phi(x) = (1/4) \psi(x) + (1/4) \psi(A^{-1}x) + (1/8) \psi(A^{-2}x);
\end{equation}
using that $d(\LG)=1/4$, and we remark that $\hat \phi$ is a piecewise
linear function since this is the case for $\hat\psi$. The conclusion
from Theorem~\ref{thm:constr-dual-wavelet-Rn-general} is that $\psi$
and $\phi$ generate dual frames $\setsmall{\dila \tran[k/2] \psi}_{j,k
  \in \Z}$ and $\setsmall{\dila \tran[k/2] \phi}_{j,k \in \Z}$ for
$L^2(\R^2)$.

The frame bounds can be found using Lemma~\ref{thm:bessel2-Rn} since the
series~(\ref{eq:14}) and (\ref{eq:lemma-suff-lower}) are finite sums on
$E$; for $\setsmall{\dila \tran[k/2] \psi}$ one finds $C_1=4/3$
and $C_2 = 4$. 
\end{example}

When the result on constructing pairs of dual wavelet frames is
written in the generality of
Theorem~\ref{thm:constr-dual-wavelet-Rn-general}, it is not always
clear how to choose the set $E$ and the lattice $\LG$. In
Example~\ref{ex:adhoc-construc-Rn} we showed how this can be done for
the quincunx dilation matrix and constructed a pair of dual frame
wavelets. In Section~\ref{sec:spec-case-Rn} and
Theorem~\ref{thm:constr-dual-wavelet-Rn-special} we specify how to
choose $E$ and $\LG$ for general dilations. The issue of exhibiting
functions $\psi$ satisfying the condition (\ref{eq:dyadic-part2-Rn})
is addressed in Section~\ref{sec:examples-rn}.

In one dimension, however, it is straightforward to make good choices
of $E$ and $\LG$ as is seen by the following corollary of
Theorem~\ref{thm:constr-dual-wavelet-Rn-general}. The corollary
unifies the construction procedures in Theorem~2 and Proposition~1
from \cite{lemvig_constr_paris} in a general procedure.
\begin{corollary}
  \label{thm:constr-dual-wavelet-gen-phi}
  Let $\dd \in \N _0$, $a>1$, and $\psi \in L^2(\R)$. Suppose that $\hat
  \psi$ is a bounded, real-valued function with $\supp \hat \psi \subset
  \itvccs{-a^{c}}{-a^{c-\dd-1}} \cup \itvccs{a^{c-\dd-1}}{a^{c}}$ for some $c
  \in \Z$, and that
  \begin{equation}
    \label{eq:dyadic-part2}
    \sum_{j \in \mathbb{Z}} \hat \psi (a^j \xi) =1 
    \quad \text{for \almoste\ } \xi \in \R.
  \end{equation}
  Let $b_j \in \C$ for $j=-\dd, \dots, \dd$, let $m = - \min\setprop{j}{\{b_j
    \neq 0\}}$, and define the function $\phi$ by
  \begin{align}
    \label{eq:dual-generator1}
    \phi(x) = \sum_{j=-m}^{\dd} b_j a^{-j} \psi(a^{-j}x) \quad
    \text{for } x \in \R.
  \end{align}
  Let $b \in \itvoc{0}{a^{-c}(1+a^m)^{-1}}$. If $b_0 = b$ and $b_j +
  b_{-j} = 2b$ for $j=1,2, \dots, \dd$, then $\psi$ and $\phi$
  generate dual frames $\setsmall{\dila[a^j] \tranl[bk] \psi}_{j,k \in
    \Z}$ and $\setsmall{\dila[a^j] \tranl[bk] \phi}_{j,k \in \Z}$ for
  $L^2(\R)$.
\end{corollary}
\begin{proof}
  In Theorem~\ref{thm:constr-dual-wavelet-Rn-general} for $n=1$ and
  $A=a$ we take $E=\itvccs{-a^c}{-a^{c-1}} \cup \itvccs{a^{c-1}}{a^c}$
  as the multiplicative tiling set for $\setprop{a^j}{j \in \Z}$. The
  assumption on the support of $\hat \psi$ becomes 
    \[  \supp \hat \psi \subset \bigcup_{j=0}^\dd a^{-j}(E) = \itvccs{-a^{c}}{-a^{c-\dd-1}}
    \cup \itvccs{a^{c-\dd-1}}{a^{c}}. \]
  Moreover, since
    \[ \bigcup_{j=0}^\dd a^{-j}(E) \subset \itvcc{-a^c}{a^c}, \quad
    \bigcup_{j=-m}^{2\dd} a^{-j}(E) \subset \itvcc{-a^{c+m}}{a^{c+m}}, \]
and
\[ (\itvcc{-a^c}{a^c}+\kstar ) \cap \itvcc{-a^{c+m}}{a^{c+m}} =
\emptyset \quad \text{for $\abs{\kstar} \ge a^c + a^{c+m} = a^c
  (1+a^m)$,}\] 
the choice $\LG^\ast = b^{-1} \Z$ for $b^{-1} \ge a^c
(1+a^m)$ satisfies equation (\ref{eq:13}). This corresponds to $\LG
=b\Z$ for $0 < b \le a^{-c} (1+a^m)^{-1}$.
\end{proof}

The assumptions in Corollary~\ref{thm:constr-dual-wavelet-gen-phi}
imply that $m \in \{0,1, \dots , \dd \}$; we note that in case $m=0$, the
corollary reduces to \cite[Theorem 2]{lemvig_constr_paris}.

\section{A special case of the construction procedure}
\label{sec:spec-case-Rn}

We aim for a more automated construction procedure than what we have
from Theorem~\ref{thm:constr-dual-wavelet-Rn-general}, in particular,
we therefore need to deal with good ways of choosing $E$ and $\LG$.
The basic idea in this automation process will be to choose $E$ as a
dilation of the difference between $I_\ast$ and $B^{-1}(I_\ast)$,
where $I_\ast$ is the unit ball in a norm in which the matrix $B=A^t$
is expanding ``in all directions''; we will make this statement
precise in Section~\ref{sec:dilation-matrix}. This idea is
instrumental in the proof of
Theorem~\ref{thm:constr-dual-wavelet-Rn-special}.

\subsection{Some results on expansive matrices}
\label{sec:dilation-matrix}
We need the following well-known equivalent conditions for a
(non-singular) matrix being expansive.
\begin{proposition}\label{thm:expansive-matrix-equiv}
  For $B \in GL_n(\R)$ the following assertions are equivalent:
  \begin{compactenum}[(i)]
    \item $B$ is expansive, \ie all eigenvalues $\lambda_i$ of $B$
      satisfy $\abs{\lambda_i}>1$. \label{item:pairsRn1}
    \item For any norm $\enorm{\,\cdot\,}$ on $\Rn$ there are
      constants $ \lambda > 1$ and $c\ge 1$ such that \[
      \enormsmall{B^jx} \ge 1/c \lambda^j \enorm{x}
      \qquad \text{for all $j \in \N_0$}, \] for any $x \in
      \Rn$. \label{enu:almost-expanding} 
    \item There is a Hermitian norm $\enorm[\ast]{\,\cdot\,}$ on $\Rn$ and a constant
    $\lambda > 1$ such that \[ 
    \enormsmall[\ast]{B^jx} \ge  \lambda^j \enorm[\ast]{x}  \qquad
    \text{for all $j \in \N_0$}, \] 
    for any $x \in \Rn$. \label{enu:really-expanding}
  \item   $\mathcal{E} \subset \lambda \mathcal{E} \subset B \mathcal{E}$ for
      some ellipsoid  $\mathcal{E}=\setpropsmall{x \in \Rn}{\enorm{Px} \le
        1}$, $P \in GL_n(\R)$, and $\lambda >1$.   \label{enu:expanding-ellipsoid}
  \end{compactenum}
\end{proposition}

By Proposition~\ref{thm:expansive-matrix-equiv} we have that for a
given expansive matrix $B$, there exists a scalar product with the
induced norm $\enorm[\ast]{\,\cdot\,}$ so that 
\[
\enorm[\ast]{Bx} \ge \lambda \enorm[\ast]{x} \quad \text{for } x\in
\Rn,
\] 
holds for some $\lambda >1$. We say that $\enorm[\ast]{\,\cdot\,}$ is
a norm associated with the expansive matrix $B$. Note that such a norm
is not unique; we will follow the construction as in the proof of
\cite[Lemma 2.2]{MR2004e:42023}, so let $c$ and $\lambda$ be as in
\enumref{enu:almost-expanding} in
Proposition~\ref{thm:expansive-matrix-equiv} for the standard
Euclidean norm with $1<\lambda < \abs{\lambda_i}$ for $i =1, \dots,
n$, where $\lambda_i$ are the eigenvalues of $B$. For $k \in \N$
satisfying $k>2 \ln c / \ln \lambda$ we introduce the symmetric,
positive definite matrix $K \in GL_\nn(\R)$:
\begin{equation}
  \label{eq:15}
  K = I + (B^{-1})^t B^{-1} + \dots +  (B^{-k})^t B^{-k}.
\end{equation}
The scalar product associated with $B$ is then defined by
$\innerprod[\ast]{x}{y}= x^t K y$. It might not be effortless to
estimate $c$ and $\lambda$ for some given $B$, but it is obvious that
we just need to pick $k \in \N$ such that $ B^t K B - \lambda^2 K $
becomes positive semi-definite for some $\lambda >1$ since this
corresponds to $\innerprod{KBx}{Bx} \ge \lambda^2 \innerprod{Kx}{x}$,
that is, $\enorm[\ast]{Bx}^2 \ge \lambda^2 \enorm[\ast]{x}^2$ for all
$x \in \Rn$.

We let $I_\ast$ denote the unit ball in the Hermitian norm  
$\enorm[\ast]{\,\cdot\,} = \enormsmall{K^{1/2}\cdot}$ associated
with $B$, \ie
\begin{equation}
 I_\ast = \setprop{x\in \Rn}{\enorm[\ast]{x} \le 1} = \setprop{x\in
  \Rn}{\enormsmall{K^{1/2}x} \le 1} =  \setprop{x\in
  \Rn}{x^t Kx \le 1},\label{eq:25}
\end{equation}
and we let $O_\ast$ denote the annulus
\[ O_\ast =I_\ast \setminus B^{-1}(I_\ast).\] 
The ringlike structure of $O_\ast$ is guaranteed by the fact that $B$
is expanding in all directions in the $\enorm[\ast]{\,\cdot\,}$ norm, \ie
\begin{equation}
 I_\ast \subset \lambda I_\ast \subset B(I_\ast), \qquad \lambda >1, \label{eq:9}
 \end{equation}
 which is \enumref{enu:expanding-ellipsoid} in Proposition
 \ref{thm:expansive-matrix-equiv}. We note that by an orthogonal
 substitution $I_\ast$ takes the form $\setpropsmall{x \in \Rn}{\mu_1
   \tilde x_1^2+\dots+\mu_\nn \tilde x_\nn^2 \le 1}$, where $\mu_i$ are
 the positive eigenvalues of $K$ and $x= Q\tilde x$ with the $i$th
 column of $Q \in
 O(\nn)$ comprising of the $i$th eigenvector of $K$. The annulus $O_\ast$ is a bounded multiplicative tiling set
 for $\setpropsmall{B^j}{j\in \Z}$. This is a consequence of the
 following result.
\begin{lemma}\label{thm:dilated-annuli}
Let $B \in GL_\nn(\R)$ be an expansive matrix. For $x\neq 0$ there is a unique $j
\in \Z$ so that $B^j x \in O_\ast$; that is, 
\begin{equation}
  \label{eq:8}
  \Rn \setminus \{0\} = \bigcup_{j \in \Z} B^j(O_\ast) \quad \text{  with disjoint union.}
\end{equation}
\end{lemma}
\begin{proof}
 From equation~(\ref{eq:9}) we know that $\{B^l(I_\ast) \}_{l\in\Z}$ is a
 nested sequence of subsets of $\Rn$, thus
 \[ B^{l}(I_\ast) \setminus B^{l-1} (I_\ast) = B^l(O_\ast), \quad l
 \in\Z ,\] 
 are disjoint sets. Since $\enorm[\ast]{B^{-j}x} \le
 \lambda^{-j}\enorm[\ast]{x}$ and $\enorm[\ast]{B^{j}x} \ge
 \lambda^{j}\enorm[\ast]{x}$ for $j \ge 0$ and $\lambda > 1$, we
 also have 
\begin{gather*}
  \bigcup_{m=-l+1}^{l} B^m(O_\ast) = B^l(I_\ast) \setminus B^{-l}(I_\ast)
 = \setprop{x\in\Rn}{\enormsmall[\ast]{B^{-l}x} \le 1 \text{ and }
   \enormsmall[\ast]{B^{l}x} > 1} 
 \\  \supset \setprop{x\in\Rn}{\lambda^{-l}\enorm[\ast]{x} \le 1 \text{ and }
   \lambda^l  \enorm[\ast]{x} > 1} =  \setprop{x\in\Rn}{\lambda^{-l} < \enorm[\ast]{x} \le \lambda^l}.
\end{gather*}
Taking the limit $l \to \infty$ we get (\ref{eq:8}).
\end{proof}

\begin{example}
\label{ex:hermitian-norm}
Let the following dilation matrix be given
\begin{equation}
  \label{eq:17}
  A =
  \begin{pmatrix}
    3 & -3 \\ 1 & 0
  \end{pmatrix}.
\end{equation}
Here we are interested in the transpose matrix $B=A^t$ with
eigenvalues $\mu_{1,2}=3/2 \pm i \sqrt{3}/2$, hence $B$ is an
expansive matrix with $\abs{\mu_{1,2}}=\sqrt{3}>1$. The dilation
matrix $B$ is not expanding in the standard norm
$\enorm[2]{\,\cdot\,}$ in $\Rn$, \ie $I_2  \not\subset B(I_2)$, as 
shown by Figure~\ref{fig:ex3-euclidean-norm}.
    \begin{figure}[ht]
      \centering
      \includegraphics[scale=.35]{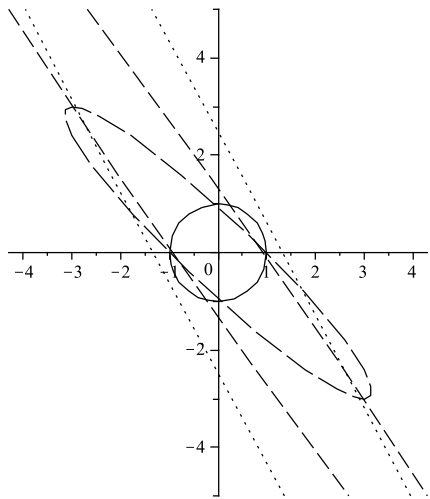}
      \caption{Boundaries of the sets $I_2$, $B(I_2)$, $B^2(I_2)$, and
        $B^3(I_2)$ marked by solid, long dashed, dashed, and dotted
        lines, respectively. Note that $I_2 \setminus B(I_2)$ is
        non-empty, and even $I_2 \setminus B^2(I_2)$ is non-empty.}
      \label{fig:ex3-euclidean-norm}
    \end{figure}
    In order to have $B$ expanding the unit ball we need to use the
    Hermitian norm from \enumref{enu:really-expanding} in
    Proposition~\ref{thm:expansive-matrix-equiv} associated with $B$.
    In \eqref{eq:15} we take $k=2$ so that the real, symmetric,
    positive definite matrix $K$ is
\begin{equation*}
  K = I + (B^{-1})^t B^{-1} +  (B^{-2})^t B^{-2} =
  \begin{pmatrix}
    28/9 & 16/9 \\ 16/9 & 8/3
  \end{pmatrix},
\end{equation*}
and let $\innerprod[\ast]{x}{y} := x^t Ky$. The choice $k=2$ suffices since
it makes $ B^t K B - \lambda^2 K $ semi-positive definite for
$\lambda =1.03$ and thus 
\[ \enorm[\ast]{Bx} \ge \lambda \enorm[\ast]{x}, \qquad x \in \R^2, \]
holds for $\lambda=1.03$. 

Figure~\ref{fig:ex3-adaptive-norm} and
\ref{fig:ex3-adaptive-norm-zoom} illustrate that $B$ indeed expands
the Hermitian norm unit ball $I_\ast$ in all directions.
    \begin{figure}[ht]
      \centering
      \includegraphics[scale=.4]{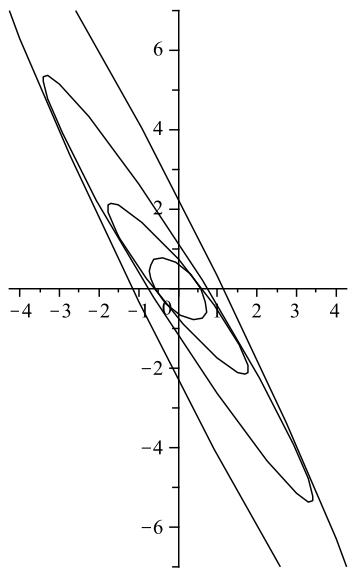}
   \caption{The unit ball $I_\ast$ in the Hermitian norm
     $\enorm[\ast]{\,\cdot\,}$ associated with $B$ and its dilations  $B(I_\ast), B^2(I_\ast),
     B^3(I_\ast)$. Only the boundaries are marked. }
      \label{fig:ex3-adaptive-norm}
    \end{figure}
We also remark that the Hermitian norm with $k=1$ will not make the dilation
matrix $B$ expanding in $\Rn$; in this case we have a situation
similar to Figure~\ref{fig:ex3-euclidean-norm}. 
    \begin{figure}[th!]
      \centering
      \includegraphics[scale=.4]{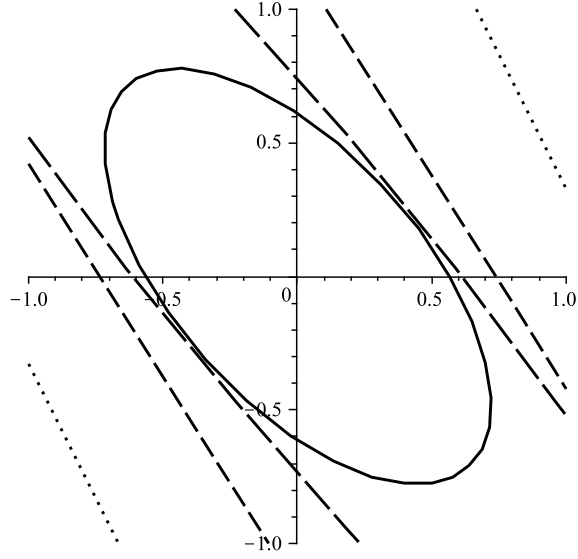}
   \caption{A zoom of Figure~\ref{fig:ex3-adaptive-norm}. Boundaries
     of the sets $I_\ast$,  $B(I_\ast)$, 
     $B^2(I_\ast)$, and $B^3(I_\ast)$ marked by solid,
     long dashed, dashed, and dotted lines, respectively.}
      \label{fig:ex3-adaptive-norm-zoom}
    \end{figure}
  \end{example}

\subsection{A crude lattice choice}
\label{sec:choos-transl-latt}
Let us consider the setup in
Theorem~\ref{thm:constr-dual-wavelet-Rn-general} with the set $E =
B^c(O_\ast)$ for some $c\in\Z$, where the norm
$\enorm[\ast]{\,\cdot\,}=\enormsmall{K^{1/2}\cdot}$ is associated with
$B$. Let $\mu$ be the smallest eigenvalue of $K$ such that $\ell =
\sqrt{1/\mu}$ is the largest semi-principal axis of the ellipsoid
$I_\ast$, \ie $\ell = \max_{x \in I_\ast}\enorm[2]{x}$. Then we can
take any lattice $\LG = P\Zn$, where $P$ is a non-singular matrix
satisfying
\begin{equation}
\normsmall[2]{P} \le \frac{1}{ \ell  \norm[2]{A^c} (1 +
  \norm[2]{A^{\underline{m}}}) }, \label{eq:16}
\end{equation}
as our translation lattice in Theorem~\ref{thm:constr-dual-wavelet-Rn-general}. To see this, recall that we are looking for
a lattice $\LG^\ast$ such that, for $\kstar \in \LG^\ast \setminus \{0\}$,
\begin{align}
  \supp \hat \phi \cap \supp \hat \psi(\cdot \pm \kstar) =
  \emptyset. \label{eq:19}
\end{align}
 For our choice of $E$ we
find that $\supp \hat \phi
\subset B^{c+\underline{m}}(I_\ast)$ and $\supp \hat \psi \subset
B^{c}(I_\ast)$. Since
\begin{align*} 
  \enorm[2]{B^{c+\underline{m}}x} \le \norm[2]{B^{c+\underline{m}}}
  \enorm[2]{x} \le  \norm[2]{B^{c+\underline{m}}} \ell \qquad \text{for any $x \in I_\ast$}, 
\end{align*}
and similar for $B^cx$, we have the situation in (\ref{eq:19}) whenever
$\enorm[2]{\kstar} \ge \ell (\norm[2]{A^c} +
\norm[2]{A^{c+\underline{m}}})$. Here we have used that
for the $2$-norm $\norm[2]{A}= \norm[2]{B}$. 
For $z \in \Zn$ we have
\begin{equation*}
 \enorm[2]{z} \le \normsmall[2]{P^t} \,\enormsmall[2]{(P^t)^{-1}z} = \normsmall[2]{P} \,\enormsmall[2]{(P^t)^{-1}z},
 \end{equation*}
therefore, by $\enorm[2]{z} \ge 1$ for $z \neq 0$,  we have 
\[ \enorm[2]{(P^t)^{-1}z} \ge \frac{1}{\norm[2]{P}} \qquad \text{for
  $z \in \Z \setminus \{0\}$.}\] Now, by assuming that $P$ satisfies
(\ref{eq:16}), we have
\begin{equation*}
\enorm[2]{\kstar} =  \enormsmall[2]{(P^t)^{-1}z} \ge 1/ \normsmall[2]{P} \ge
\ell  \norm[2]{A^c} (1 +  \norm[2]{A^{\underline{m}} }) \ge l (\norm[2]{A^c} +
\norm[2]{A^{c+\underline{m}}})
\end{equation*}
for $0 \neq \kstar = (P^t)^{-1}z \in (P^t)^{-1}\Zn=\LG^\ast$,
hence the claim follows.

A lattice choice based on (\ref{eq:16}) can be rather crude, and
produces consequently a wavelet system with unnecessarily many
translates. From equation~(\ref{eq:16}) it is obvious that any lattice
$\LG=P \Zn$ with $\norm{P}$ sufficiently small will work as
translation lattice for our pair of generators $\psi$ and $\phi$.
Hence, the challenging part is to find a sparse translation lattice
whereby we understand a lattice $\LG$ with large determinant
$d(\LG):=\abs{\dtm{P}}$. In the dual lattice system this corresponds
to a dense lattice $\LG^\ast$ with small volume $d(\LG^\ast)$ of the
fundamental parallelotope $I_{\LG^\ast}$ since $d(\LG)d(\LG^\ast)=1$.
In Theorem~\ref{thm:constr-dual-wavelet-Rn-special} in the next
section we make a better choice of the translation lattice compared to
what we have from (\ref{eq:16}).

Using a crude lattice  approach as above,  we can easily
transform the translation lattice to the integer lattice if we allow
multiple generators. We pick a matrix $P$ that satisfies condition
(\ref{eq:16}) and whose inverse is integer valued, \ie  
 $Q:=P^{-1} \in GL_n(\Z)$.  The conclusion from
 Theorem~\ref{thm:constr-dual-wavelet-Rn-general} is that
 $\setsmall{\dila \tran[ Q^{-1}k] \psi}_{j 
  \in \Z,  k \in \Zn}$ and $\setsmall{\dila \tran[Q^{-1}k] \phi}_{j \in
  \Z, k \in \Zn}$ are dual frames. 
The order of the quotient group $Q^{-1}\Zn/\Zn$ is $\abs{\dtm{Q}}$, so  
let  $\setpropsmall{d_i }{i=1, \dots, \abs{\dtm{Q}}}$ denote a
complete set of representatives of the quotient group, and define \[
\Psi = \setprop{\tran[d_i] \psi}{i = 
   1, \dots, \abs{\dtm Q}},  \quad  \Phi = \setprop{\tran[d_i] \phi}{i =
   1, \dots, \abs{\dtm Q}}.\] 
Since $\setsmall{\dila \tran[ Q^{-1}k] \psi}_{j 
  \in \Z,  k \in \Zn} = \setsmall{\dila \tran[ k] \psi}_{j 
  \in \Z,  k \in \Zn, \psi \in \Psi}$ and likewise for the dual frame,
the statement follows.


\subsection{A concrete version of
  Theorem~\protect{\ref{thm:constr-dual-wavelet-Rn-general}} }
\label{sec:special-case-Rn}
We list some standing assumptions and conventions for this section.

\paragraph{\bf General setup.} We assume $A \in GL_\nn(\R)$ is 
expansive. Let $\enorm[\ast]{\,\cdot\,} =
\innerprod[\ast]{\,\cdot\,}{\,\cdot\,}^{1/2}$ be a Hermitian norm as
in \enumref{enu:really-expanding} in
Proposition~\ref{thm:expansive-matrix-equiv} associated with $B =
A^t$, let $I_\ast$ denote the unit ball in the $\enorm[\ast]{\,\cdot\,}$-norm,
and let $K \in GL_\nn(\R)$ be the symmetric, positive definite matrix
such that $\innerprod[\ast]{x}{y} = y^t Kx$. Let $\Lambda :=
\mathrm{diag}(\lambda_1, \dots, \lambda_\nn)$, where $\{\lambda_i\}$
are the eigenvalues of $K$, and let $Q \in O(\nn)$ be such that the
spectral decomposition of $K$ is $Q^t K Q = \Lambda$.

The following result is a special case of
Theorem~\ref{thm:constr-dual-wavelet-Rn-general}, where we, in
particular, specify how to choose the translation lattice $\LG$. Since
we in Theorem~\ref{thm:constr-dual-wavelet-Rn-special} define $\LG$,
it allows for a more automated construction procedure.

\begin{theorem}\label{thm:constr-dual-wavelet-Rn-special} 
  Let $A, I_\ast, K, Q, \Lambda$ be as in the general setup. Let $\dd \in \N
  _0$ and $\psi \in L^2(\Rn)$. Suppose that $\hat \psi$ is a bounded,
  real-valued function with $\supp \hat \psi \subset B^c (I_\ast)
  \setminus B^{c-d-1} (I_\ast)$ for some $c \in \Z$, and that
  \eqref{eq:dyadic-part2-Rn} holds. Take $\LG = (1/2) A^c Q
  \sqrt{\Lambda} \Zn$. Then the function $\psi$ and the function
  $\phi$ defined by
\begin{align}
  \label{eq:dual-generator2-Rn}
  \phi(x) = d(\LG) \left [ \psi(x) + 2 \sum_{j=0}^{\dd}
    \abs{\dtm{A}}^{-j} \psi(A^{-j}x)\right ] \quad
  \text{for } x \in \Rn,
\end{align}
 generate dual frames $\setsmall{\dila \tran \psi}_{j \in \Z, \gamma \in
  \LG}$ and $\setsmall{\dila \tran \phi}_{j \in \Z, \gamma \in \LG}$ for
$L^2(\Rn)$
  \end{theorem}

  \begin{remark}
    Note that $d(\LG) = 2^{-n} \abs{\dtm{A}}^c (\lambda_1 \cdots
    \lambda_\nn)^{1/2} $ and $\sqrt{\Lambda} =
    \mathrm{diag}(\sqrt{\lambda_1}, \dots, \sqrt{\lambda_\nn})$.
  \end{remark}

\begin{proof}
  The annulus $O_\ast$ is a bounded multiplicative tiling set for the dilations
  $\setprop{B^j}{j\in\Z}$ by Lemma~\ref{thm:dilated-annuli}, hence
  this is also the case for $B^c(O_\ast)$ for $c\in \Z$. The support
  of $\hat\psi$ is $\supp \hat \psi \subset B^c (I_\ast) \setminus
  B^{c-d-1} (I_\ast) = \cup_{j=0}^\dd B^{c-j}(O_\ast)$. Therefore we
  can apply Theorem~\ref{thm:constr-dual-wavelet-Rn-general} with
  $E=B^c(O_\ast)$, $b_j=2$ and $b_{-j}=0$ for $j =1, \dots, \dd$ so
  that $\underline{m}=0$ and $\overline{m}=\dd$. The only thing left
  to justify is the choice of the translation lattice
  $\LG$. 
  We need to show that condition~(\ref{eq:13}) with $\underline{m}=0$
  and $\overline{m}=\dd$ in
  Theorem~\ref{thm:constr-dual-wavelet-Rn-general} is satisfied by
  $\LG^\ast = 2 B^{c} Q \Lambda^{-1/2} \Zn$. By the orthogonal
  substitution $x = Q \tilde x$ the quadratic form $x^t K x$ of
  equation (\ref{eq:25}) reduces to
\[ \lambda_1 \tilde x_1^2 + \dots +  \lambda_\nn \tilde x_\nn^2, \]
where $\lambda_i>0$, hence in the $\tilde x = Q^t x$ coordinates $I_\ast$ is
given by 
\[ \tilde I_\ast = \setprop{\tilde x \in \Rn}{\biggl(\frac{\tilde
    x_1}{1/\sqrt{\lambda_1}}\biggr)^2+ \dots + \biggl(\frac{\tilde
    x_\nn}{1/\sqrt{\lambda_\nn}}\biggr)^2 < 1}\]
   which is an ellipsoid with semi axes $\tfrac{1}{\sqrt{\lambda_1}},
   \dots, \tfrac{1}{\sqrt{\lambda_\nn}}$.  
   Therefore, in the $\tilde x$ coordinates,  
\[ (\tilde I_\ast + \kstar) \cap \tilde I_\ast = \emptyset  \qquad \text{for
  $0 \neq \kstar \in 2 \Lambda^{-1/2} \Zn $},\]
or, in the $x$ coordinates, 
\begin{equation*} 
  ( I_\ast + \kstar) \cap  I_\ast = \emptyset  \qquad \text{for
  $0 \neq \kstar \in 2 Q \Lambda^{-1/2} \Zn $.}
\end{equation*}
By applying $B^c$ to this relation it becomes
\begin{equation} \label{eq:20}
  \bigl( B^c(I_\ast) + \kstar \bigr) \cap  B^c(I_\ast) = \emptyset  \qquad \text{for
  $0 \neq \kstar \in \LG^\ast = 2 B^c Q \Lambda^{-1/2} \Zn $,}
\end{equation}
whereby we see that condition~(\ref{eq:13}) is satisfied with
$\underline{m}=0$ and $\LG^\ast = 2 B^{c} Q \Lambda^{-1/2} \Zn$. The
dual lattice of $\LG^\ast$ is $\LG = 1/2 A^{-c} Q \Lambda^{1/2} \Zn$.
It follows from Theorem~\ref{thm:constr-dual-wavelet-Rn-general} that
$\psi$ and $\phi$ generate dual frames for this choice of the translation
lattice.
\end{proof}

The frame bounds for the pair of dual frames $\setsmall{\dila \tran
  \psi}_{j \in \Z, \gamma \in \LG}$ and $\setsmall{\dila \tran
  \phi}_{j \in \Z, \gamma \in \LG}$ in
Theorem~\ref{thm:constr-dual-wavelet-Rn-special} can be given
explicitly as
\begin{align*}
  C_1 &= \frac{1}{d(\LG)}\inf_{\xi \in B^{c-\dd}(O_\ast)} \sum_{j=0}^{\dd}
  \left(\hat \psi (B^{j} \xi)\right)^2, & &C_2 =
  \frac{1}{d(\LG)}\sup_{\xi \in B^{c-\dd}(O_\ast)} \sum_{j=0}^{\dd} \left(\hat
    \psi (B^{j} \xi)\right)^2, 
  \intertext{and} 
  C_1 &= \frac{1}{d(\LG)}\inf_{\xi \in B^{c-\dd}(O_\ast)}
  \sum_{j=-\dd}^{\dd} \left(\hat \phi (B^{j}
  \xi)\right)^2, & &C_2 = \frac{1}{d(\LG)}\sup_{\xi \in B^{c-\dd}(O_\ast)}
  \sum_{j=-\dd}^{\dd} \left(\hat \phi (B^{j}
  \xi)\right)^2 ,
\end{align*}
respectively. The frame bounds do not depend on the specific structure
of $\LG$, but only on the determinant of $\LG$; in particular, the
condition number $C_2/C_1$ is independent of $\LG$.

To verify these frame bounds, we note that equation (\ref{eq:20})
together with the fact $\supp \hat \psi$, $\supp \hat \phi \subset
B^c(I_\ast)$ imply that
\begin{align*}
  \hat \psi(\xi) \hat \psi(\xi + \kstar) = \hat \phi(\xi) \hat
  \phi(\xi + \kstar) = 0 \qquad \text{for \almoste\  $\xi \in \Rn$
    and $\kstar \in \LG^\ast \setminus \{0\}$}.
\end{align*}
Therefore, by equations~(\ref{eq:18}) and (\ref{eq:28}) with
$E=B^c(O_\ast)$, $\underline{m}=0$ and $\overline{m}=\dd$, we have
\begin{align*}
  \sum_{j\in\Z} \sum_{\kstar \in \LG^\ast } \abs{ \hat \psi(B^j \xi)
    \hat \psi(B^j \xi + \kstar)} =  \sum_{j\in\Z} \abs{ \hat \psi(B^j \xi)}^2
  =   \sum_{j=0}^{\dd} \left( \hat \psi(B^j \xi) \right)^2,
\end{align*}
and 
\begin{align*}
  \sum_{j\in\Z} \sum_{\kstar \in \LG^\ast } \abs{ \hat \phi(B^j \xi)
    \hat \phi(B^j \xi + \kstar)} =  \sum_{j\in\Z} \abs{ \hat \phi(B^j \xi)}^2
  =   \sum_{j=-\dd}^{\dd} \left( \hat \phi(B^j \xi) \right)^2,
\end{align*}
for $\xi \in B^{c-\dd}(O_\ast)$. The stated frame bounds follow from
Lemma~\ref{thm:bessel2-Rn}.

\begin{example}
Let $A$ and $K$ be as in Example~\ref{ex:hermitian-norm}. 
The eigenvalues of $K$ are $\lambda_1  = (26 + 2\sqrt{65})/9 \approx
4.7$ and $\lambda_2  = (26 - 2\sqrt{65})/9 \approx 1.1$.
Let the normalized (in the standard norm) eigenvectors of $K$ be columns of
$Q \in O(2)$ and $\Lambda = \mathrm{diag}(\lambda_1, \lambda_2)$,
hence $Q^t K Q = \Lambda$.  
By the orthogonal transformation 
  $x = Q \tilde x$ the Hermitian norm unit ball $I_\ast$ becomes
\[ \tilde I_\ast = \setprop{\tilde x \in \R^2}{\biggl(\frac{\tilde
    x_1}{1/\sqrt{\lambda_1}}\biggr)^2+ \biggl(\frac{\tilde
    x_2}{1/\sqrt{\lambda_2}}\biggr)^2 < 1} \subset I_2 \]
 which is an ellipse with semimajor axis $1/ \sqrt{\lambda_2} \approx 0.95$
 and semiminor axis $1/ \sqrt{\lambda_1} \approx 0.46$. 
Since 
$ \Lambda^{-1/2} = \mathrm{diag}(1/\sqrt{\lambda_1},
1/\sqrt{\lambda_2})$, we have
\[ \abs{(\tilde I_\ast + \kstar) \cap \tilde I_\ast} = 0  \qquad \text{for
  $0 \neq  \kstar \in 2 \Lambda^{-1/2} \Z^2 $}.\]
By the orthogonal  substitution back to $x$ coordinates, we get
\[ \abs{( I_\ast + \kstar) \cap  I_\ast} = 0  \qquad \text{for
  $0 \neq \kstar \in 2 Q \Lambda^{-1/2} \Z^2 $}.\]

Suppose that $\hat \psi$ is a bounded, real-valued
  function with $\supp \hat \psi \subset B^c (I_\ast) \setminus
  B^{c-\dd-1} (I_\ast)$ for $c = 1$ that satisfies the $B$-dilative partition~\eqref{eq:dyadic-part2-Rn}.
Since $c=1$ we need to take $\LG^\ast = 2 B^{1} Q
\Lambda^{-1/2} \Z^2$ and $\LG = 1/2 A^{-1} Q \Lambda^{1/2} \Z^2$, see
Figure~\ref{fig:ex3-gamma-star-c1} and \ref{fig:ex3-gamma-c1}.
    \begin{figure}[ht]
      \centering
      \includegraphics[scale=.35]{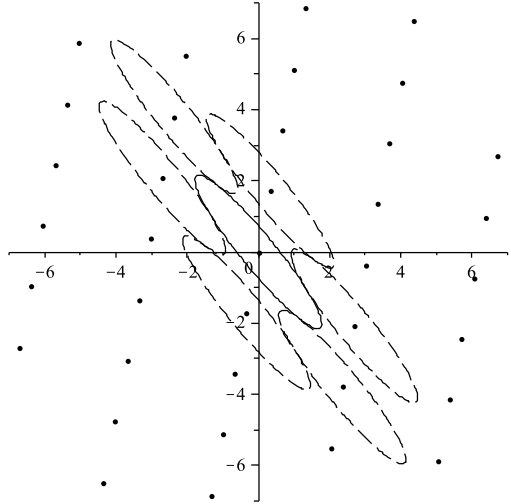}
      \caption{The dual lattice $\LG^\ast = 2B^{c}Q
        {\LL}^{-1/2} \Z^2$ for $c=1$ is shown by dots, and the boundary of the
        set $B^c(I_\ast)$ by a solid line. Boundaries of the set
        $B^c(I_\ast)$ translated to several different $\kstar \in
        \LG^\ast \setminus \{0\}$ are shown with dashed lines. Recall that
        $\supp \hat \psi, \supp \hat \phi \subset B^c(I_\ast)$, hence
        $\supp \hat \phi \cap \supp \hat \psi(\cdot + \kstar) =
        \emptyset$ for $\kstar \in
        \LG^\ast \setminus \{0\}$. }
      \label{fig:ex3-gamma-star-c1}
    \end{figure}
    \begin{figure}[ht]
      \centering
      \includegraphics[scale=.35]{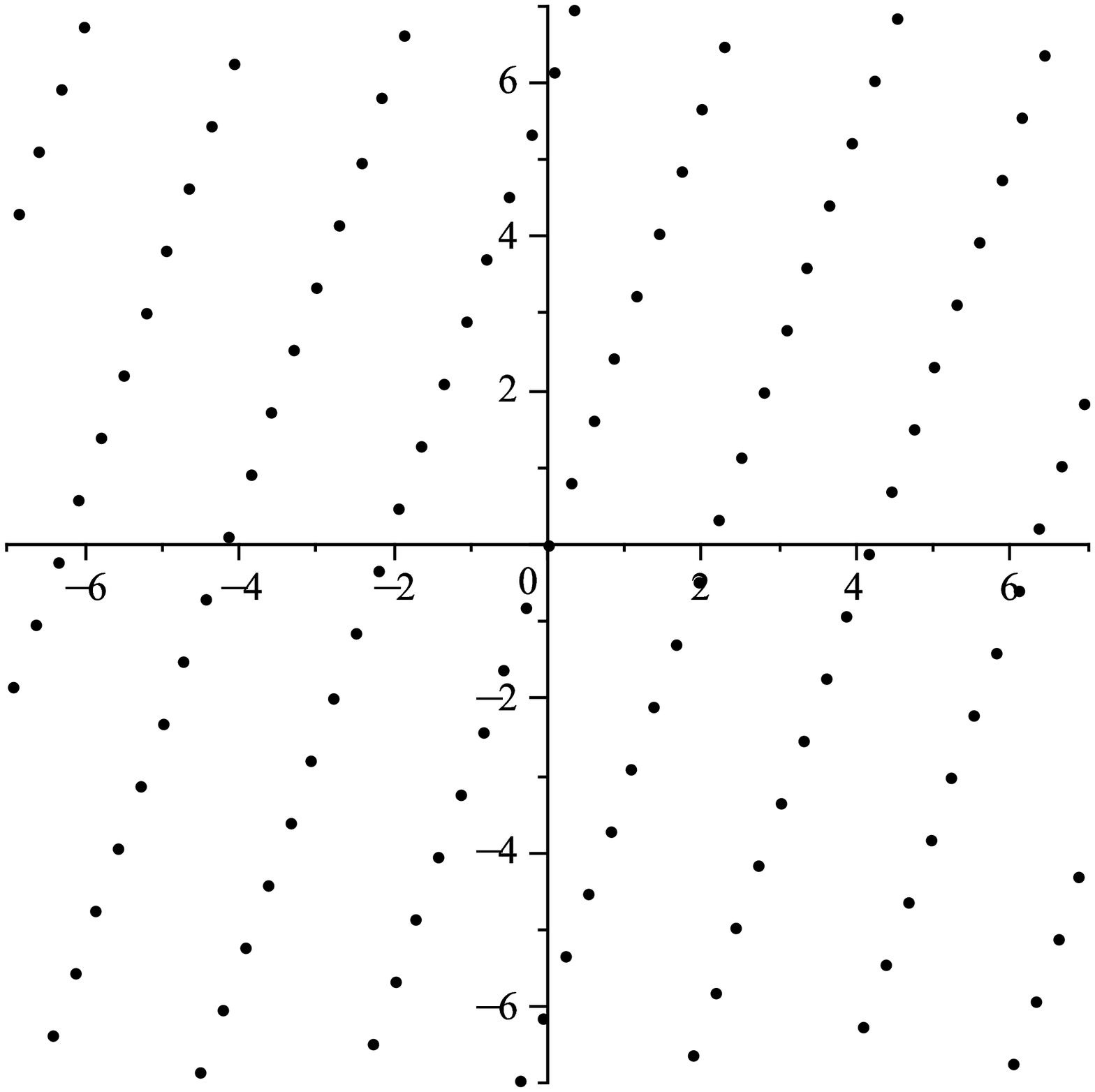}
   \caption{The translation lattice $\LG = (1/2) A^{c}Q {\LL}^{1/2}
     \Z^2$ for $c=1$.}
      \label{fig:ex3-gamma-c1}
    \end{figure}
\end{example}

\subsection{An alternative lattice choice}
\label{sec:optim-latt-choice}
Let the setup up and assumptions be as in
Theorem~\ref{thm:constr-dual-wavelet-Rn-special}, except for the
lattice $\LG$ which we want to choose differently. As in
Section~\ref{sec:choos-transl-latt} the dual lattice $\LG^\ast$ needs
to satisfy (\ref{eq:19}) for $\kstar \in \LG^\ast \setminus
\{0\}$.  We want to choose $\LG^\ast$ as dense as possible since this will
make the translation lattice $\LG$ as sparse as possible and the
wavelet system with as few translates as possible. Since $\supp \hat
\psi, \supp \hat \phi \subset B^c(I_\ast)$, we are looking for
lattices $\LG^\ast$ that packs the ellipsoids $B^c(I_\ast) +
\kstar$, $\kstar \in \LG^\ast$, in a non-overlapping,
optimal way. By the coordinate transformation $\hat x =
\Lambda^{-1/2}Q^t B^{-c}x$, the ellipsoid $B^c(I_\ast)$ turns into the
standard unit ball $I_2$ in $\Rn$. This calculations are as follows.
\begin{align*}
  B^c(I_\ast) &= \setprop{B^cx}{\enorm[\ast]{x}^2 \le 1} =
  \setprop{x}{\enormsmall[2]{K^{1/2}B^{-c}x}^2 \le 1} \\
  &= \setprop{x}{\enorm[2]{K^{1/2}B^{-c}B^c Q \Lambda^{-1/2}\hat x}^2
    \le 1} \\
  &= \setprop{x}{\einnerprod[2]{\hat x}{\Lambda^{-1/2} Q^t K
       Q\Lambda^{-1/2} \hat x} \le 1} =\setprop{x}{\enorm[2]{\hat x}^2
     \le 1},
\end{align*}
and we arrive at a standard sphere packing problem with lattice
arrangement of non-overlapping unit $n$-balls. The proportion of the
Euclidean space $\Rn$ filled by the balls is called the density of the
arrangement, and it is this density we want as high as possible.  

Taking $\LG$ as in Theorem~\ref{thm:constr-dual-wavelet-Rn-special}
corresponds to a square packing of the unit $n$-balls $I_2 + k$ by the
lattice $2 \Zn$, \ie $k \in 2\Zn$. The density of this packing is $V_n
2^{-n}$, where $V_n$ is the volume of the $n$-ball: $V_{2n}=\pi^n/(n!)$ and
$V_{2n+1}=(2^{2n+1}n!\pi^n)/(2n+1)!$.
This is not the densest packing of balls in $\Rn$ since there exists
a lattice with density bigger than $1.68 n2^{-n}$ for each $n\neq
1$ \cite{0030.34602}; a slight improvement of this lower bound was
obtained in \cite{0776.52006} for $n>5$. Moreover, the densest lattice
packing of hyperspheres is known up to dimension $8$, see
\cite{MR1648172}; it is precisely this dense lattice we want to use in
place of $2 \Zn$ (at least whenever $n \le 8$).

In $\R^2$ Lagrange proved that the hexagonal packing, where each ball
touches $6$ other balls in a hexagonal lattice, has the highest
density $\pi / \sqrt{12}$. Hence using $P \Z^2$ with
\[ P =
\begin{pmatrix}
  2 & 0 \\ 1 & \sqrt{3}
\end{pmatrix}
\]
instead of $2 \Z^2$ improves the packing by a factor of
\[ \frac{\pi/\sqrt{12}}{\pi/2^2} = 4/\sqrt{12} = 2/\sqrt{3}.\] It is
easily seen that this factor equals the relation between the area of
the fundamental parallelogram of the two lattices $\abs{\dtm{2I_{2
      \times 2}}}/\abs{\dtm{P}}$. In
Figure~\ref{fig:ex3-gamma-star-c1} we see that each ellipse only
touches $4$ other ellipses corresponding to the square packing
$2\Zn$; in the optimal packing each ellipse touch $6$ others. In
$\R^3$ Gauss proved that the highest density is $\pi / \sqrt{18}$
obtained by the hexagonal close and face-centered cubic packing; here
each ball touches $12$ other balls.

\section{Dilative partition of unity}
\label{sec:examples-rn}
With Theorem~\ref{thm:constr-dual-wavelet-Rn-special} at hand the only
issue left is to specify how to construct functions satisfying the
partition of unity (\ref{eq:dyadic-part2-Rn}) for any given expansive
matrix. In the two examples of this section we outline possible ways
of achieving this.

\subsection{Constructing a partition of unity}
\label{sec:constr-part-unity}
As usual we fix the dimension $n \in \N$ and the expansive matrix $B \in
GL_\nn(\R)$. 
  In the examples in this section we construct functions satisfying the assumptions
  in Theorem~\ref{thm:constr-dual-wavelet-Rn-special}, that is, a real-valued
  function $g \in L^2(\Rn)$ with $\supp g \subset B^c(I_\ast)
  \setminus B^{c-\dd-1}(I_\ast)$ for some $c \in \Z$ and $\dd \in \N_0$ so that
the $B$-dilative partition 
\begin{equation}
     \sum_{j \in \mathbb{Z}} g (B^j \xi) =1 
\quad \text{for \almoste\ }\xi \in \Rn,
 \label{eq:21}
 \end{equation}
holds. 

In the construction we will use that the radial coordinate of the
surface of the ellipsoid $\partial B^j(I_\ast)$, $j \in \Z$, can be
parametrized by the $\nn-1$ angular coordinates $\theta_1, \dots,
\theta_{\nn-1}$.
  The radial coordinate expression will  be of the form
  $h(\theta_1, \dots, \theta_{\nn-1})^{-1/2}$ for some positive,
  trigonometric function $h$, where $h$ is bounded away from zero and
  infinity with the specific form of $h$ depending on the dimension
  $n$ and the length and orientation of the ellipsoid axes.

  We illustrate this with the following example in $\R^4$. We want to
  find the radial coordinate $r$ of the ellipsoid \[\setprop{x \in
    \R^4}{(x_1/\ell_1)^2 + (x_2/\ell_2)^2 +(x_3/\ell_3)^2 + (x_4/\ell_4)^2 =1},
  \qquad \ell_i > 0,\, i=1,2,3,4, \] as a function the angular coordinates
  $\theta_1, \theta_2$ and $\theta_3$. We express $x=(x_1, x_2, x_3,
  x_4) \in \R^4$ in the hyperspherical coordinates $(r,\theta_1,
  \theta_2, \theta_3) \in \{0\} \cup \R_+ \times [0,\pi] \times
  [0,\pi] \times [0, 2\pi)$ as follows:
\begin{align*}
  x_1 &= r \cos \theta_1, && x_2 = r \sin \theta_1 \cos \theta_2, \\ 
  x_3 &= r \sin \theta_1 \sin \theta_2 \cos \theta_3, && x_4 = r  \sin
  \theta_1 \sin \theta_2  \sin \theta_3.  
\end{align*}
Then we substitute $x_i$, $i=1,\dots,4$, in the expression above and factor
out $r^2$ to obtain $r^2 f(\theta_1, \theta_2, \theta_{3}) =1$, where
  \begin{align}\label{eq:27}
    f(\theta_1,\theta_2, \theta_3) &= \ell_1^{-2}\cos^2 \theta_1 +  \ell_2^{-2}\sin^2 \theta_1
    \cos^2 \theta_2 \\ &\phantom{=}\,+  \ell_3^{-2}\sin^2 \theta_1 \sin^2 \theta_2 \cos^2 \theta_3 +  \ell_4^{-2} \sin^2
  \theta_1 \sin^2 \theta_2  \sin^2 \theta_3. \nonumber
  \end{align}
The conclusion is that $r=r(\theta_1,\theta_2, \theta_3) = f(\theta_1, \theta_2,
  \theta_{3})^{-1/2}$.

\begin{example}\label{exa:d-one}
  For $\dd=1$ in Theorem~\ref{thm:constr-dual-wavelet-Rn-special} we
  want $g \in C^s_0(\Rn)$ for any given $s \in \N \cup \{0\}$. The
  choice $\dd=1$ will fix the ``size'' of the support of $g$ so that $\supp
  g \subset B^c(I_\ast) \setminus B^{c-2}(I_\ast)$ for some $c\in \Z$.
  Now let $r_1=r_1(\theta_1, \dots, \theta_{\nn-1})$ and
  $r_2=r_2(\theta_1, \dots, \theta_{\nn-1})$ denote the radial
  coordinates of the surface of the ellipsoids $\partial
  B^{c-1}(I_\ast)$ and $\partial B^c(I_\ast)$ parametrized by $\nn-1$
  angular coordinates $\theta_1, \dots, \theta_{\nn-1}$, respectively.

  Let $f$ be a continuous function on the annulus
  $S=\overline{B^c(O_\ast)}$ satisfying $f \vert _{\partial
    B^{c-1}(I_\ast)}=1$ and $f \vert _{\partial
    B^{c}(I_\ast)}=0$.
  Using the parametrizations $r_1, r_2$ of the surfaces of the two
  ellipsoids and fixing the $\nn-1$ angular coordinates we realize
  that we only have to find a continuous function $f:\itvcc{r_1}{r_2}
  \to \R$ of one variable (the radial coordinate) satisfying
  $f(r_1)=1$ and $f(r_2)=0$. For example the general function $f \in
  C^0(S)$ of $d$ variables can be any of the functions below:
\begin{subequations}
 \label{eq:def-of-f-radius}
   \begin{align}
    f(x)&= f(r,\theta_1, \dots, \theta_{\nn-1}) =
    \frac{r_2-r}{r_2-r_1},  \label{eq:22}\\
    f(x)&= f(r,\theta_1, \dots, \theta_{\nn-1}) = \frac{(r_2-r)^2}{(r_2-r_1)^3} (2 (r -r_1)+r_2-r_1),\label{eq:23}  \\ 
    f(x)&=  f(r,\theta_1, \dots, \theta_{\nn-1}) = \tfrac{1}{2}+\tfrac{1}{2}\cos{\pi (\tfrac{r-r_1}{r_2-r_1})} \label{eq:24},
\end{align}
\end{subequations}
where $r = \enorm{x} \in \itvcc{r_1}{r_2}$, $\theta_1, \dots,
\theta_{\nn-2} \in \itvcc{0}{\pi}$, and $\theta_{\nn-1} \in
\itvco{0}{2\pi}$; recall that $r_1=r_1(\theta_1, \dots, \theta_{\nn-1})$
and $r_2=r_2(\theta_1, \dots, \theta_{\nn-1})$. In definitions
\eqref{eq:23} and \eqref{eq:24} the function $f$ even 
belongs to $C^1(S)$.

Define $g \in L^2(\R)$ by:
\begin{equation}
  \label{eq:def-C1-example-Rn}
  g(x) =
  \begin{cases}
    1-f(Bx)  \quad &\text{for } x \in B^{c-1}(I_\ast)
  \setminus B^{c-2}(I_\ast), \\ 
   f(x)  \quad &\text{for } x \in B^{c}(I_\ast)
  \setminus B^{c-1}(I_\ast), \\
   0   \quad & \text{otherwise.} 
  \end{cases}
\end{equation}
This way $g$ becomes a $B$-dilative partition of unity with $\supp
g \subset  B^{c}(I_\ast)
  \setminus B^{c-2}(I_\ast) $, so we can
apply Theorem~\ref{thm:constr-dual-wavelet-Rn-special} with $\hat
\psi=g$ and $d=2$. 
\end{example}

We can simplify the expressions for the radial coordinates $r_1, r_2$
of the surface of the ellipsoids $\partial B^{c-1}(I_\ast)$ and
$\partial B^c(I_\ast)$ from the previous example by a suitable
coordinate change. The idea is to transform the ellipsoid $
B^{c-1}(I_\ast)$ to the standard unit ball $I_2$ by a first coordinate
change $\tilde x = \Lambda^{1/2} Q^t B^{-c+1} x$. This will transform
the outer ellipsoid $B^c(I_\ast)$ to another ellipsoid. A second and
orthogonal coordinate transform $\hat x = Q_\prime^t \tilde x$ will
make the semiaxes of this new ellipsoid parallel to the coordinate
axes, leaving the standard unit ball $I_2$ unchanged. Here $Q_\prime$
comes from the spectral decomposition of $A^{-1}B^{-1}$, \ie
$A^{-1}B^{-1} = Q_\prime^t \Lambda_\prime Q_\prime$. In the $\hat x$
coordinates $r_1=1$ is a constant and $r_2=f^{-1/2}$ with $f$ of the
form (\ref{eq:27}) for $n=4$ and likewise for $n \neq
4$. 

In the construction in Example~\ref{exa:d-one} we assumed that $\dd=1$.  The next
example works for all $\dd \in \N$; moreover, the constructed function
will belong to $C^\infty_0(\Rn)$.   
\begin{example}\label{exa:all-d}
For sufficiently small $\delta >0$ define $\Delta_1, \Delta_2 \subset \Rn$ by
\begin{align*}
  \Delta_1 &= B^{c-\dd-1}(I_\ast) + \unitball{0}{\delta}, \\
  \Delta_2 &+  \unitball{0}{\delta} = B^{c}(I_\ast). 
\end{align*}
This makes $\Delta_2 \setminus \Delta_1$ a subset of the
annulus $B^{c}(I_\ast) \setminus B^{c-\dd-1}(I_\ast)$; it is exactly the
subset, where points less than $\delta$ in
distance from the boundary have been removed, or in other words
\[  \Delta_2 \setminus \Delta_1 + \unitball{0}{\delta} =
B^{c}(I_\ast) \setminus B^{c-\dd-1}(I_\ast). 
\]
For this to hold, we of course need to take $\delta >0$ sufficiently
small, \eg such that $\Delta_1 \subset r
\Delta_1 \subset \Delta_2$ holds for some $r>1$. 

Let $h\in C^\infty _0(\Rn)$ satisfy $\supp h = \unitball{0}{1}$, $h \ge
0$, and $\int h \,\mathrm{d}\mu =1$, and define $h_\delta = \delta^{-d}
h(\delta^{-1}\cdot)$. By convoluting the characteristic function on
$\Delta_2 \setminus \Delta_1$ with $h_\delta$ we obtain a smooth
function living on the annulus $ B^{c}(I_\ast) \setminus
B^{c-\dd-1}(I_\ast)$. So let $p \in C^\infty _0(\Rn)$ be defined by 
\[ p = h_\delta \ast \chi_{\Delta_2 \setminus \Delta_1},
\]
and note that  $\supp p = B^{c}(I_\ast) \setminus B^{c-\dd-1}(I_\ast)$ since 
$\supp h_\delta = \unitball{0}{\delta}$. Normalizing the function $p$
in a proper way will give us the function $g$ we are looking for. We
will normalize $p$ by the function $w$: 
\begin{equation*}
  w(x) = \sum_{j \in \mathbb{Z}} p (B^j x). 
\end{equation*}
For a fixed $x\in \Rn \setminus \{0\}$ this sum has either $d$ or $d+1$ nonzero
terms, and $w$ is therefore bounded away from
$0$ and $\infty$:
\begin{equation*}
  \exists c, C > 0 \, \colon c < w(x) < C \quad \text{for all } x\in
  \Rn \setminus \{0\},
\end{equation*}
hence we can define a function $g\in C^\infty_0(\Rn)$ by
\begin{equation}
  g(x) = \frac{p(x)}{w(x)} \quad \text{for } x\in \Rn \setminus \{0\},
  \quad \text{and,} \quad g(0)=0.
  \label{eq:definition-g-Rn}
  \end{equation}
The function $g$ will be an almost everywhere $B$-dilative partition of unity as is seen by
using the $B$-dilative periodicity of $w$: 
\begin{align*}
  \sum_{j \in \mathbb{Z}} g (B^j x) = \sum_{j \in \mathbb{Z}}
  \frac{p (B^j x)}{w(B^j x)}
  = \sum_{j \in \mathbb{Z}}
  \frac{p(B^j x)}{w(x)} 
= \frac{1}{w(x)} \sum_{j \in \mathbb{Z}}p (B^j x) = 1.
\end{align*}
Since $p$ is supported on the annulus $B^{c}(I_\ast) \setminus
B^{c-\dd-1}(I_\ast)$, we can simplify the definition in
(\ref{eq:definition-g-Rn}) to get rid of the infinite sum in the
denominator; this gives us the following expression
\[  g(x) = p(x) / \sum_{j=-\dd}^{\dd}p(B^j x) \qquad \text{for
} x \in\Rn \setminus \{0\}. \]

We can obtain a more explicit expression for $p$ by the following
approach. Let $r_1=r_1(\theta_1, \dots, \theta_{\nn-1})$ and
$r_2=r_2(\theta_1, \dots, \theta_{\nn-1})$ denote the radial coordinates
of the surface of the ellipsoids $\partial B^{c-\dd-1}(I_\ast)$ and
$\partial B^c(I_\ast)$ parametrized by $\nn-1$ angular coordinates
$\theta_1, \dots, \theta_{\nn-1}$, respectively. 
Finally, let $p \in C^\infty_0(\Rn)$ be defined by
\[ p(x) = \eta(\enorm{x}-r_1)\,\eta(r_2-\enorm{x}), \quad \text{with }
r_1= r_1(\theta_1, \dots, \theta_{\nn-1}) \text{ and } r_2=
r_2(\theta_1, \dots, \theta_{\nn-1})\] 
where $\theta_1, \dots,
\theta_{\nn-1}$ can be found from $x$, and
\begin{equation*}
    \eta(x)= \begin{cases} \mathrm{e}^{-1/x} \quad  &x>0, \\ 0 \quad  &x\le 0.
    \end{cases}
  \end{equation*}
\end{example}

\end{document}